\documentclass[a4paper,10pt]{amsart}

\usepackage{cite}
\usepackage{amsmath,amsthm}
\usepackage{amssymb,latexsym}
\usepackage{amsfonts}
\usepackage{hyperref}
\usepackage[ngerman,english]{babel}

\newcommand{\N}{\mathbb{N}}
\newcommand{\K}{\mathbb{K}}
\newcommand{\from}{\colon}
\DeclareMathOperator{\Id}{Id}

\newtheorem{theorem}{Theorem}

\title{Graph Laplacians do not generate strongly continuous semigroups}
\author{T. Kalmes, C. Schumacher}
\thanks{\noindent 2010 \textit{Mathematics Subject Classification.} 47D06, 05C63, 46A04\\
\textit{Key words and phrases}: graph Laplacians, strongly continuous semigroup on locally convex spaces}

\begin{document}

\begin{abstract}
We show that for graph Laplacians $\Delta_G$ on a connected locally finite simplicial undirected graph~$G$ with countable infinite vertex set~$V$
none of the operators $\alpha\,\Id+\beta\Delta_G, \alpha,\beta\in\K,\beta\ne0$,
generate a strongly continuous semigroup on~$\K^V$ when the latter is equipped with the product topology.
\end{abstract}

\maketitle

\section{Introduction}

Let $G=(V,E)$ be a connected locally finite simplicial undirected graph with a countably infinite vertex set~$V$ and the set of edges $E\subseteq V\times V$.
Two vertices $v,w\in V$ of~$G$ are adjacent, $v\sim w,$ if $(v,w)\in E$, i.\,e.\ $(v,w)$ is an edge of~$G$.
Recall that a graph is locally finite if for every vertex~$v$ of~$G$ the number $\operatorname{deg}(v)$ of adjacent vertices is finite and that~$G$ is connected if for any pair of different vertices $v,w\in V$ there is a path connecting~$v$ and~$w$, i.e.\ a finite number of edges $(v_0,v_1),\ldots,(v_{n-1},v_n)$ of~$G$ with $v\in\{v_0,v_n\}$ and $w\in\{v_0,v_n\}$.
The length of a path is the number of its edges.
The graph~$G$ is simplicial if it does have no multiple edges, which we excluded implicitly by $E\subseteq V\times V$, and no loops loops, i.\,e.\ $(v,v)$ is not an edge of~$G$ for any vertex~$v$.
Finally,~$G$ is undirected if $(w,v)$ is an edge of~$G$ whenever $(v,w)$ is an edge of~$G$.

A Laplacian on~$G$ is a linear operator acting on the real valued functions on~$V$, $\Delta_G\from\K^V\to\K^V$, defined by
\[\forall\,v\in V\colon\,\Delta_G\,f(v)=\sum_{v\sim w}\gamma_{v,w}(f(v)-f(w))\]
with positive weights $\gamma_{v,w}>0$.
Common choices are $\gamma_{v,w}=1$ or $\gamma_{v,w}=1/\operatorname{deg}(v)$.
If we equip~$\K^V$ with the product topology it is clear that~$\K^V$ is a locally convex vector space on which~$\Delta_G$ is a continuous linear mapping.
$\Delta_G$ has been investigated e.g.\ in \cite{CeCo09}, \cite{CeCoDo12}, and \cite{Kalmes14}.

The definitions and most basic results for semigroups on locally convex spaces~$X$ are the same as for Banach spaces, we refer to \cite{Komatsu64}, \cite{Komura68}, \cite{Ouchi73}, \cite{Yosida65}.
A strongly continuous semigroup~$T$ on~$X$ is thus a morphism from the semigroup $([0,\infty),+)$ to that of continuous linear operators $(L(X),\circ)$ such that all orbits $t\mapsto T_t(x)$, $x\in X$, are continuous.

Although under rather weak assumptions on~$X$ many results from the Banach space setting carry over to strongly continuous semigroups of locally convex spaces there are some crucial differences, e.\,g.\ in general not every continuous linear operator on~$X$ is the generator of a strongly continuous semigroup on~$X$.

The purpose of this short note is to proof the following theorem.
For the definition of a quasiadjacency matrix with respect to $G$ see section 2.

\begin{theorem}
For every $\lambda\in\K$ and each quasiadjacency matrix~$B$ with respect to~$G$ the continuous linear operator $\lambda\Id-B$ does not generate a strongly continuous semigroup on~$\omega$.
In particular $\alpha\Id+\beta\Delta_G$ is not the generator of a strongly continuous semigroup on $\K^V$ for any $\alpha,\beta\in\K$, $\beta\ne0$.
\end{theorem}

\section{Proof of Theorem 1}

Enumeration of the vertices $V=\{v_k;\,k\in\N\}$ of~$G$ with $v_k\ne v_j$ for $j\ne k\in\N$ gives an isomorphism of~$\K^V$ onto~$\K^\N$.
In order to keep notation simple, for a linear mapping $A\from\K^V\to\K^V$ we denote by~$A$ also the linear operator on~$\K^\N$ induced by the isomporphism between~$\K^V$ and~$\K^\N$.

We equip~$\K^\N$ with its usual Fr\'echet space topology, i.\,e.\ the locally convex topology defined by the increasing fundamental system of seminorms $(p_k)_{k\in\N}$ given by
\[\forall\,f=(f_j)_{j\in\N}\in\K^\N\colon\,p_k(f)=\sum_{j=1}^k\lvert f_j\rvert.\]
As usual, we denote this Fr\'echet space by~$\omega$.
$\omega$ is the projective limit of the projective spectrum of Banach spaces $(\K^n)_{n\in\N}$ with surjective linking maps $\pi_m^n\from\K^m\to\K^n,\pi_m^n(x_1,\ldots,x_m)=(x_1,\ldots,x_n), m\ge n$, thus $\omega$ is a quojection. Moreover, let $\pi_m\from\K^\N\to \K^m,\pi_m((x_n)_{n\in\N})=(x_1,\ldots,x_m)$. 

Using that $G$ is connected and locally finite, a straightforward calculation shows that $A\from\omega\to\omega$ is continuous if and only if the inducing $A\from\K^V\to \K^V$ has finite hopping range, i.\,e.\ if for every $v\in V$ there is $n\in\N$ such that
\[\forall\,f,g\in\K^V\colon f_{|U_n(v)}=g_{|U_n(v)}\implies A(f)(v)=A(g)(v),\]
where $U_n(v)$ is the $n$-neighbourhood of $v$, i.\,e.\ $U_n(v)$ is the union of $\{v\}$ with the set of all endpoints of paths in~$G$ starting in~$v$ and of length not exceeding~$n$.
Obviously,~$\Delta_G$ has finite hopping range and is thus is a continuous linear mapping on~$\omega$.

Finally, we call a matrix $B=(b_{k,l})\in\K^{\N\times\N}$ \textit{quasiadjacency matrix with respect to $G$} if there is $\alpha\in\K$ such that $\alpha\,b_{k,l}\ge 0$ for all $k,l$ and $b_{k,l}\ne 0$ precisely when $v_k$ and $v_l$ are adjacent.
Obviously, $\Delta_G=\Id-B$ for some quasiadjacency matrix~$B$ with respect to~$G$.
Moreover, for each quasiadjacency matrix~$B$ with respect to~$G$ and for every $\lambda\in\K$ the linear operator $\lambda\,\Id-B$ has finite hopping range and is thus continuous on~$\omega$.\\  

Since $\omega$ is a quojection, according to \cite[Theorem 2]{FrJoKaWe14} $\lambda\,\Id-B$ generates a strongly continuous semigroup if and only if
\begin{equation}\label{characterisation}
\forall\,n\in\N\,\exists\,m\in\N\,\forall\,k\in\N_0, x\in\omega:\,\big(\pi_m(x)=0\Rightarrow\pi_n((\lambda\,\Id-B)^k x)=0\big).
\end{equation}
Denoting the unit vectors by $e_l=(\delta_{k,l})_{k\in\N}$ we have for $l\ge 2$ and $k\in\N$
\[\pi_1\big((\lambda\,\Id-B)^k e_l\big)=\sum_{j=1}^k \binom{k}{j}(-1)^j\lambda^{k-j}\pi_1(B^j e_l).\]
Since $B$ is a quasiadjacency matrix with respect to $G$ it follows that $\pi_1(B^j e_l)$ does not vanish precisely when there is a path of length $j$ in $G$ connecting $v_1$ and $v_l$.

Now fix $m$ and $l>m$ so that $\pi_m(e_l)=0$. Moreover, set
\[k:=\min\{\text{length $j$ of a path connecting $e_1$ and $e_l$}\}\text.\]
Since $G$ is connected, $k\in\N$ and
\[\pi_1\big((\lambda\,\Id-B)^k e_l\big)=\sum_{j=1}^k \binom{k}{j}(-1)^j\lambda^{k-j}\pi_1(B^j e_l)=(-1)^k\pi_1(B^k e_l)\ne 0.\]
Since $m$ is arbitrary, (\ref{characterisation}) is not fulfilled for $n=1$ proving the theorem.\hfill$\square$

\begin{small}
\selectlanguage{ngerman}
\noindent{\textsc{Technische Universit\"at Chemnitz,
Fakult\"at f\"ur Mathematik, 09107 Chemnitz, Germany}}
\par
\noindent{\textit{E-mail addresses: thomas.kalmes@mathematik.tu-chemnitz.de,\\ christoph.schumacher@mathematik.tu-chemnitz.de}}
\end{small}
\end{document}